\documentclass[12pt]{amsart}

\usepackage{amsfonts,amsmath,latexsym,amssymb,verbatim,amsbsy}
\usepackage{amsthm}

\usepackage{pstricks}

\numberwithin{equation}{section}

\theoremstyle{plain}
\newtheorem{THEOREM}{Theorem}[section]

\newtheorem{theorem}[THEOREM]{Theorem}

\newtheorem{lemma}[THEOREM]{Lemma}

\theoremstyle{definition}

\newtheorem{definition}[THEOREM]{Definition}

\theoremstyle{remark}

\newcommand{\thm}[1]{Theorem~\ref{#1}}
\newcommand{\lem}[1]{Lemma~\ref{#1}}

\newcommand{\N}{\ensuremath{\mathbb{N}}}   
\newcommand{\Z}{\ensuremath{\mathbb{Z}}}   
\newcommand{\R}{\ensuremath{\mathbb{R}}}   
\newcommand{\C}{\ensuremath{\mathbb{C}}}   
\newcommand{\T}{\ensuremath{\mathbb{T}}}   

\def \a {\alpha}
\def \b {\beta}
\def \d {\delta}
\def \g {\gamma}
\def \e {\varepsilon}
\def \f {\varphi}

\def \k {\kappa}
\def \l {\lambda}

\def \i {\iota}
\def \n {\nabla}
\def \s {\sigma}

\def \Th {\Theta}

\def \ba {{\bf a}}
\def \bA {{\bf A}}
\def \bB {{\bf B}}

\def \bF {{\bf F}}
\def \bP {{\bf P}}
\def \bPi {\boldsymbol{\Pi}}

\def \bs {\boldsymbol{\sigma}}
\def \bI {{\bf I}}
\def \id {{\bf id}}
\def \br {{\bf r}}

\def \bS {{\bf S}}

\def \bH {{\bf H}}

\def \bL {{\bf L}}
\def \bG {{\bf G}}
\def \bT {{\bf T}}
\def \bS {{\bf S}}

\def \bU {{\bf U}}
\def \bV {{\bf V}}

\def \bFi {{\bf \Phi}}

\def \bs {\boldsymbol{\sigma}}
\def \bt {\boldsymbol{\tau}}
\def \bl {\boldsymbol{\lambda}}

\def \torus {{\T^n}}
\def \rn {{\dot{\R}^n}}
\def \zn {{\dot{\Z}^n}}

\def \encon {{L}_\mathcal{F}^2}

\def \frbundle {\mathcal{F}}

\def \p {\partial}
\def \ra {\rightarrow}
\def \ss {\subset}

\newcommand{\der}[2]{(#1 \cdot \nabla) #2}

\newcommand{\ress}[1]{r_{\mathrm{ess}}(#1)}
\newcommand{\sym}[1]{\mathcal{S}^{#1}}

\newcommand{\tA}[2]{\tilde{A}_{#1,\, #2}}
\newcommand{\tB}[2]{\tilde{B}_{#1,\, #2}}

\DeclareMathOperator{\supp}{supp} %
\DeclareMathOperator{\op}{Op} %
\DeclareMathOperator{\re}{Re} %
\DeclareMathOperator{\im}{Im} %
\DeclareMathOperator{\Ker}{Ker} %

\begin{document}

\title[The unstable spectrum of the {N}avier-{S}tokes]{The unstable spectrum of the {N}avier-{S}tokes operator in the limit of vanishing viscosity}
\author{Roman Shvydkoy and Susan Friedlander}
\address{Department of Mathematics, Statistics, and Computer Science, University of Illinois at Chicago,
(M/C 249), Chicago, IL 60607} %
\email{susan@math.uic.edu; shvydkoy@math.uic.edu}
\thanks{Friedlander's research is partially supported by NSF grants DMS-0202767 and
DMS-0503768. The authors thank Mathematisches Forschungsinstitut
Oberwolfach for hospitality during our work.}
\date{\today}%

\maketitle

\begin{abstract}
A general class of linear advective PDEs, whose leading order term
is of viscous dissipative type, is considered. It is proved that
beyond the limit of the essential spectrum of the underlying
inviscid operator, the eigenvalues of the viscous operator, in the
limit of vanishing viscosity, converge precisely to those of the
inviscid operator. The general class of PDEs includes the equations
of incompressible fluid dynamics. Hence eigenvalues of the
Navier-Stokes operator converge in the inviscid limit to the
eigenvalues of the Euler operator beyond the essential spectrum.
\end{abstract}

\section{Introduction}

The equations of motion governing an incompressible fluid with
viscosity $\e$ are the Navier-Stokes equations
\begin{subequations}\label{NSE}
\begin{align}
\frac{\p q_\e}{\p t} &= - \der{q_\e}{q_\e} - \n p_\e + \e \Delta
q_\e
+ F_\e, \\
\n \cdot q_\e &= 0,
\end{align}
\end{subequations}
where $q_\e$ denotes the $n$-dimensional velocity vector, $p_\e$
denotes the pressure and $F_\e$ is an external force vector. Here
$n$ can be any integer with $n\geq 2$, but the case $n = 3$ is of
the most interest.

The same equations with zero viscosity are the Euler equations
\begin{subequations}\label{EE}
\begin{align}
\frac{\p q}{\p t} &= - \der{q}{q} - \n p, \\
\n \cdot q &= 0.
\end{align}
\end{subequations}

An important connection between the Euler and the Navier-Stokes
systems is the behavior of \eqref{NSE} in the limit of vanishing
viscosity (i.e. $\e \ra 0$). This limit is likely to be crucial in
the understanding of many physical problems of fluid flow, such as
the transition to turbulence. It is clear, since the types of the
two systems are very different (\eqref{NSE} is parabolic and
\eqref{EE} is degenerate hyperbolic), that the limit of vanishing
viscosity is a subtle and singular limit. There are a number of
partial results for the nonlinear system as $\e \ra 0$. The history
of such results is briefly surveyed in the appendix of the book of
Temam \cite{Temam}.

Open questions remain even for the linearized problem. In this
present paper we address the connections between the spectra of the
linearized Navier-Stokes operators in the inviscid limit and the
spectrum of the linearized Euler operator. The results are closely
tied to issues of linear, and even nonlinear, instabilities for
fluid flows (c.f. Yudovich \cite{Yudovich}).

Let $u(x,\e)$ be an arbitrary steady solution of \eqref{NSE}:
\begin{subequations}\label{ss}
\begin{align}
0 &= - \der{u}{u} - \n P + \e \Delta u + F_\e, \\
\n \cdot u &= 0.
\end{align}
\end{subequations}
We assume that $u(x,\e)$ and $F_\e$ are infinitely smooth vector
valued functions on the torus $\torus$ with regular dependence on
$\e \in [0,\e_0)$ and that $\lim_{\e \ra 0} F_\e = 0$. For the sake
of simplicity we will present the proof of the theorems only for the
case where $u(x)$ has no dependence on $\e$. The more general
results follow from similar arguments.

The linearized Navier-Stokes equations for the evolution of a small
perturbation velocity $v(x,t)$ are
\begin{subequations}\label{LNSE}
\begin{align}
\frac{\p v}{\p t} &= - \der{u}{v} - \der{v}{u} - \n p + \e \Delta v, \\
\n \cdot v &= 0.
\end{align}
\end{subequations}
The corresponding linearized Euler equations are
\begin{subequations}\label{LEE}
\begin{align}
\frac{\p v}{\p t} &= - \der{u}{v} - \der{v}{u} - \n p , \\
\n \cdot v &= 0.
\end{align}
\end{subequations}

We will study general classes of differential operators on $\torus$
which include the operators of the fluid equations defined by
\eqref{LNSE} and \eqref{LEE}. We will investigate the relationship
between the unstable point spectrum of the inviscid operator and the
eigenvalues of the viscous operator in the limit of vanishing
viscosity.

For a general equilibrium $u(x)$ the Euler operator defined in
\eqref{LEE} is non self-adjoint, non elliptic and degenerate. Hence,
contrary to the case of the elliptic Navier-Stokes operator given by
\eqref{LNSE}, standard spectral results for elliptic operators do
not apply to the Euler operator. However in the past decade
considerable progress has been made in understanding the structure
of the spectrum of the Euler operator using techniques of geometric
optics. A survey of these results is given in Friedlander and
Lipton-Lifschitz \cite{FriedLipton}. In particular, Vishik
\cite{V96} obtained an explicit, and often computable, expression
for the essential spectral radius of the Euler evolution operator in
terms of a geometric quantity that can be considered as a "fluid"
Lyapunov exponent. Recently Shvydkoy \cite{Shv2005} has extended
these results to a general class of advective PDEs with
pseudodifferential bounded perturbation. The approach in which the
evolution operator is partitioned into high frequency and low
frequency parts will also be used in Section \ref{S:high} of this
present paper. We make use of the decomposition of the inviscid
operator proved in \cite{Shv2005} to obtain the analogous
decomposition of the viscous operator (see \thm{T:decomp}). This
result requires certain explicit estimates on the symbols of PDOs on
the torus that we state in the appendix.

In Section \ref{proof} we prove the main result. We first prove a
result for spectral convergence in the inviscid limit for the
semigroups. We then prove that beyond the limit of the essential
spectrum of the inviscid operator the eigenvalues of the viscous
operator converge precisely to those of the inviscid operator. A key
step is to use the decomposition established in Section \ref{S:high}
to split off a finite dimensional subspace corresponding to growing
modes. An analogous argument was used by Lyashenko and Friedlander
\cite{LyashFried} to obtain a sufficient condition for instability
in the limit of vanishing viscosity for a class of operators
satisfying certain compactness and accretive properties. The
properties required in \cite{LyashFried} do not hold in general for
the Euler operator \eqref{LEE} (although they do hold for the
coupled rotating fluid/body system as noted in \cite{LyashFried}).
The goal of this present paper is to adapt the arguments of
\cite{LyashFried} to a wider class of operators that include the
generic Euler and Navier-Stokes operators themselves.

\section{Formulation of the result}

We consider the following class of differential operators on
$\torus$:
\begin{equation}\label{le}
    \bL^\e f = -\der{u}{f} + \bA f + \e \Delta f.
\end{equation}
Here $u \in C^\infty(\torus)$ is a divergence-free time independent
vector field, $f$ takes values in $\C^d$, $d \in \N$, and $\bA$ is a
global pseudodifferential operator (PDO) on $\torus$ given by
\begin{equation}\label{}
\bA f(x) = \op[\ba] f (x)= \sum_{k \in \zn} e^{i k\cdot x} \ba(x,k)
\hat{f}(k),
\end{equation}
where $\zn = \Z^n \backslash \{0\}$. We assume that $\ba \in
\sym{0}$ (see the appendix) is a $d \times d$-matrix valued symbol,
which allows decomposition
$$
\ba = \ba_0 + \ba_1,
$$
where $\ba_0$ is $0$-homogenous in $\xi$, and $\ba_1 \in \sym{-1}$.
For instance, the linearized Navier-Stokes equation \eqref{LNSE} has
the right hand side of the form \eqref{le}, where $\bA$ has
principal symbol
\begin{equation}\label{ps}
    \ba_0(x,\xi) = \left( 2 \frac{\xi \otimes \xi}{|\xi|^2} - \id
\right)\p u(x)
\end{equation}
(see \cite{Shv2005} for derivation).

We further consider a smooth linear bundle $\frbundle$ over $\rn
=\R^n \backslash \{0\}$. We assume that $\frbundle$ is
$0$-homogenous. A function $f \in L^2(\torus) = L^2$ is said to
satisfy the frequency constraints determined by $\frbundle$ if
$\hat{f}(k) \in F(k)$, for all $k \in \Z^n$, where $F(k)$ denotes
the fiber over $k$ (we consider the fiber $F(0)$ separately). The
space of all such functions is denoted $\encon$.

Let $\{\bG_t^\e\}_{t\geq 0}$ be a $C_0$-semigroup generated by
$\bL^\e$ over $L^2$. We assume that $\bG_t^\e$ leaves $\encon$
invariant so that the equation
\begin{equation}\label{pde}
    f_t = \bL^\e f
\end{equation}
is well-posed on $\encon$.

The first order advective operator $\bL^0$ was treated in
\cite{Shv2005}. It was shown there (and previously in \cite{V96} for
the Euler equation) that the action of $\bG_t^0$ on shortwave
localized envelopes of the form
$$
f_\d(x) = b_0(x) e^{i \xi_0 \cdot x /\d}, \quad \d \ll 1
$$
is described by the asymptotic formula
\begin{equation}\label{asym}
\bG_t^0 f_\d (x) = \bB_t(\f_{-t}(x), \xi_0)f_\d(\f_{-t}(x)) + O(\d),
\end{equation}
as $\d \ra 0$. In this formula $\f_t$ is the integral flow of the
steady field $u$, and $\bB_t$ is the fundamental matrix solution of
the amplitude equation
\begin{equation}\label{b}
    b_t = \ba_0(\chi_t(x_0,\xi_0)) b,
\end{equation}
over the phase flow $\chi_t$ on $\Th = \torus \times \rn$ determined
by the bicharacteristic system
\begin{equation}\label{bich}
        \begin{cases}
          x_t = u(x) \\
          \xi_t = - \p u^\top(x) \xi
        \end{cases}
\end{equation}
subject to the frequency constraint $b_0 \in F(\xi_0)$. One can
modify the symbol $\ba_0$ in such a way that the action of $\bL^0$
on functions from $\encon$ is the same, while \eqref{b} leaves
$\frbundle$ invariant, i.e. $b(t) \in F(\xi(t))$ (see
\cite{Shv2005}). Thus one can consider \eqref{b} as a dynamical
system over the bundle $\frbundle$.

It was proved that the exponential instabilities of the amplitude
equation \eqref{b} not only cause exponential instability of the
semigroup $\bG^0$ via \eqref{asym}, but also create the essential
spectrum of the semigroup operator $\bG_t^0$ in the unstable region.
More precisely, the following formula for the essential spectral
radius holds:
\begin{equation}\label{misha}
\ress{\bG_t^0} = e^{t\mu},
\end{equation}
where $\mu$ is the maximal Lyapunov exponent of the dynamical system
\eqref{b} (see \cite{CL}). The main result of this present article
states that beyond this limit of the essential spectrum the
eigenvalues of $\bL^\e$ converge precisely to the eigenvalues of
$\bL^0$ (and, of course, by spectral mapping the same is true for
the semigroups). Even stronger, we show convergence of the
corresponding spectral subspaces.

For a closed operator $\bL$ we use the following notation:
$$
\s_a^+(\bL) = \{ \l \in \s(\bL): \re \l > a\},
$$
and we denote by $m_a(\l,\bL)$ the algebraic multiplicity of $\l$.

\begin{theorem}\label{main}
Suppose that $\s_\mu^+(\bL^0) \neq \emptyset$. Then
\begin{enumerate}
    \item[(i)] there exists $\e_0>0$ such that $\s_\mu^+(\bL^\e) \neq
    \emptyset$ for all $0\leq \e < \e_0$,
    \item[(ii)] for any $\l \in \s_\mu^+(\bL^0)$ and any sufficiently
    small $r>0$ there is $\e_r>0$ such that for all $\e < \e_r$ one
    has
    \begin{equation}\label{mult}
        m_a(\l,\bL^0) = \sum_{\substack{ \l' \in \s_\mu^+(\bL^\e) \\
        |\l - \l'| < r}} m_a(\l',\bL^\e),
    \end{equation}
    \item[(iii)] we have the limit
    \begin{equation}\label{limproj}
        \lim_{\e \ra 0} \sum_{\substack{ \l' \in \s_\mu^+(\bL^\e) \\
        |\l - \l'| < r}} \bP^\e_{\l'} = \bP_\l^0,
    \end{equation}
     where $\bP_\l^\e$ denotes the Riesz projection onto the
     spectral subspace corresponding to $\l$.
\end{enumerate}
\end{theorem}

We note that (iii) $\Rightarrow$ (ii) $\Rightarrow$ (i). So, it
suffices to prove only part (iii). The proof heavily relies on the
results of the next section, and will be finished in Section
\ref{proof}. In the appendix we state some of the general facts on
PDO's in the way that is convenient to use in the subsequent
arguments.

\section{High frequency decomposition}\label{S:high}

In this section we identify the high frequency part of the semigroup
operator $\bG_t^\e$. It is a PDO shifted by the flow $\f_t$, while
the rest is a sum of an operator of order $O(\sqrt{\e})$ and a
compact operator that behaves like a PDO of order $-1$ uniformly in
$\e$. We introduce the following notation. As before, we let $\f_t$
denote the flow generated by $u$ on $\torus$, and $\chi_t$ the phase
flow of the bicharacteristic equations \eqref{bich} on $\torus
\times \rn$. The fundamental matrix solution of the amplitude
equation \eqref{b}, which we denoted $\bB_t(x,\xi)$, is a smooth
linear cocycle over the flow $\chi_t$. We call it the $b$-cocycle.
Clearly, the $b$-cocycle is a $0$-homogenous in $\xi$ symbol of
class $\sym{0}$.

We consider the operator of composition with the inverse flow
$\f_{-t}$:
\begin{equation}\label{ }
\bFi_t f = f \circ \f_{-t},
\end{equation}
and the orthogonal projector:
\begin{equation}\label{ }
\bPi : L^2 \ra \encon,
\end{equation}
which, as one can easily see, is a Fourier multiplier with the
symbol given by the orthogonal frequency projector onto the fiber
$F(\xi)$.

The following decomposition was proved in \cite{Shv2005} in the case
of $\e =0$:
\begin{equation}\label{e:decompzero}
    \bG_t^0 = \bH_t^0 + \bU_t^0,
\end{equation}
where
\begin{equation}\label{ }
\bH_t^0 = \bPi \bFi_t \op[\bB_t]
\end{equation}
and $\bU_t^0$ is a compact operator, which behaves like a PDO of
order $-1$ (hence the asymptotic formula \eqref{asym}). For any
positive $\e$ formula \eqref{e:decompzero} can be generalized as
follows.

\begin{theorem}\label{T:decomp}
For any $0\leq t<T$ and $0 \leq \e < \e_0$ the following
decomposition holds:
\begin{equation}\label{E:decomp}
    \bG_t^\e = \bH_t^\e + \sqrt{\e}\, \bT_t^\e + \bU_t^\e,
\end{equation}
where
\begin{align}
\bH_t^\e & = \bPi \bFi_t \op[\bt_t^\e] , \\
\bt_t^\e(x,\xi) & = \bB_t(x,\xi) \exp\left\{-\e \int_0^t | \p
\f_s^{-\top}(x)\xi|^2 \, ds \right\},
\end{align}
the family $\{\bT_t^\e\}_{0 \leq \e < \e_0,\, 0\leq t< T}$ is
uniformly bounded, and $\{ \bU_t^\e \}_{0\leq \e < \e_0,\, 0\leq t<
T}$ is uniformly compact.
\end{theorem}

By a uniform compact family we mean the following.

\begin{definition}\label{D:uc} Let $\psi_N(\xi)$ be the characteristic function of the ball
$\{|\xi| < N\}$. Define the projection multiplier $\bP_N f = (\psi_N
\hat{f})^\vee$. We say that a family of operators $\{
\bU_\iota\}_{\iota \in I}$ on $L^2$, or its subspace invariant with
respect to $\bP_N$, is uniformly compact if
\begin{equation}\label{ }
\lim_{N \ra \infty} \sup_{\iota \in I} \| \bU_\iota - \bU_\iota
\bP_N\| = 0.
\end{equation}
\end{definition}

The rest of the section is devoted to the proof of \thm{T:decomp}.

First, we notice that the theorem easily reduces to the case when
$\bPi = \bI$. Indeed, consider the semigroup $\bG_t^\e$ defined on
the whole $L^2$. If \eqref{E:decomp} holds on all $L^2$, then by
applying $\bPi$ and restricting to $\encon$, we see that
\eqref{E:decomp} holds on $\encon$ too.

Using the fact that the $b$-cocycle solves the amplitude equations
\eqref{b} we find the evolution equation for $\bH_t^\e$ by
straightforward differentiation:
\begin{multline}
\frac{d}{dt}\bH_t^\e  = -\der{u}{\bH_t^\e } + \bFi_t \op[(\ba_0\circ
\chi_t) \bt_t^\e] - \label{evol}\\ - \e \bFi_t \op[|\p
\f_t^{-\top}(x)\xi|^2 \bt_t^\e ].
\end{multline}

We compare the second term on the right hand side with $\bA
\bH_t^\e$. First, the change of variables rule implies
$$
\bA \bFi_t = \bFi_t \bA'
$$
where $\bA' = \op[\ba_0 \circ \chi_t ] + \op[\ba_t']$ with $\ba_t'
\in \sym{-1}$ uniformly in $t<T$. The latter follows from the fact
that $\f_t \in C^\infty(\torus)$ uniformly in $-T < t < T$. Hence we
obtain
\begin{equation}\label{ah}
\bA \bH_t^\e = \bA \bFi_t \bS_t^\e = \bFi_t \op[\ba_0 \circ \chi_t
]\op[\bt_t^\e] +  \bFi_t \op[\ba_t'] \op[\bt_t^\e].
\end{equation}

Let us show that the symbols $\bt_t^\e$ satisfy a uniformity
condition in the $x$-variable.

\begin{lemma}\label{L:unif}
For any multi-index $\a$ there exists a constant $B_\a$ independent
of $0\leq \e<\e_0$ and $t<T$ such that
\begin{equation}\label{E:unif}
   \sup_{x\in\torus,\, \xi \neq 0} |\p_x^\a \bt_t^\e(x,\xi) | \leq
   B_\a.
\end{equation}
\end{lemma}
\begin{proof}
By the Leibnitz rule,
$$
\p_x^\a \bt_t^\e = \sum_{\a' \leq \a} \p^{\a-\a'}_x \bB_t(x,\xi)
\p_x^{\a'} \exp\left\{-\e \int_0^t | \p \f_s^{-\top}(x)\xi|^2 \, ds
\right\}.
$$
Hence, by the uniform boundedness of the $b$-cocycle,
$$
|\p_x^\a \bt_t^\e| \leq C_{t,\a} \sup_{\a' \leq \a,\, x,\, \xi}
\left|\p_x^{\a'} \exp\left\{-\e \int_0^t | \p \f_s^{-\top}(x)\xi|^2
\, ds \right\}\right|.
$$

One can check by induction that if $g = g(x_1,\ldots,x_n)$ is a
smooth function, then
$$
\p_x^\a(\exp(g)) = \exp(g) \sum (\p_x^{\g_1}g )^{l_1} \ldots
(\p_x^{\g_r}g)^{l_r}
$$
where the sum is taken over a subset of indexes satisfying
$$
|\g_1| l_1 + \ldots + |\g_r| l_r = |\a|.
$$

In our case $g = -\e \int_0^t | \p \f_s^{-\top}(x)\xi|^2 \, ds$.
Then,
$$
|(\p_x^{\g_1}g)^{l_1} \ldots (\p_x^{\g_r}g)^{l_r}| \leq C_{t,\a}
\e^{l_1+\ldots+l_r} |\xi|^{2(l_1+\ldots+l_r)}.
$$
Using that $| \p \f_t^{-\top}(x)\xi| \geq c_{t}|\xi|$, we obtain
$$
\left|\p_x^{\a'} \exp\left\{-\e \int_0^t | \p \f_s^{-\top}(x)\xi|^2
\, ds \right\}\right| \leq C_{t,\a} e^{-c_t\e|\xi|^2} \sum (\e
|\xi|^2)^{l_1+\ldots+l_r} \leq C'_{t,\a}
$$
uniformly in $\e$.
\end{proof}

Using Lemmas \ref{L:unif} and \ref{L:uc2} we immediately conclude
that the family $\bU_{t,\e}^{(1)} = \bFi_t \op[\ba_t']
\op[\bt_t^\e]$ is uniformly compact in $0\leq t<T$ and $0\leq \e
<\e_0$. By Theorem \ref{compose}, with $m_1 = m_2 = 0$ and $N = 4$,
$$
\op[\ba_0 \circ \chi_t ]\op[\bt_t^\e] = \op[\bl_t^\e],
$$
where
$$
\bl_t^\e =(\ba_0 \circ \chi_t) \bt_t^\e +\sum_{1\leq |\g| < 4}
\frac{(-1)^{|\g|}}{\g!}( \p_\xi^\g \ba_0 \circ \chi_t)(\p_x^\g
\bt_t^\e) + \br_4^{t,\e}.
$$
From the estimate on the remainder \eqref{rem} and the $\bt_t^\e$
given in \eqref{E:unif}, we see that the families of symbols
$\{\br_4^{t,\e}\}$ and $\{( \p_\xi^\g \ba_0 \circ \chi_t)(\p_x^\g
\bt_t^\e)\}$, with $|\g|\geq 1$, satisfy the assumption of Lemma
\ref{L:uc1}. Hence, they contribute a uniformly compact family
$\{\bU_{t,\e}^{(2)}\}$.

Summarizing the above, we have shown the identity
\begin{equation}\label{ah2}
    \bA \bH_t^\e = \bFi_t \op[(\ba_0 \circ \chi_t) \bt_t^\e] +
    \bU_{t,\e}^{(3)},
\end{equation}
where $\{\bU_{t,\e}^{(3)}\}$ is uniformly compact.

Let us now consider the term
\begin{equation}\label{nextterm}
- \e \bFi_t \op[|\p \f_t^{-\top}(x)\xi|^2 \bt_t^\e ]
\end{equation}
and compare it to
\begin{equation}\label{nextcomp}
    \e \Delta \bH_t^\e = -\e \op[|\xi|^2] \bFi_t \op[\bt_t^\e].
\end{equation}
By the change of variables rule, we obtain
$$
\op[|\xi|^2] \bFi_t = \bFi_t \op[|\p \f_t^{-\top}(x)\xi|^2].
$$
By \thm{compose} with $m_1 = 2$, $m_2 = 0$, $N = 6$
$$
\e \op[|\p \f_t^{-\top}(x)\xi|^2] \op[\bt_t^\e] = \op[\bl_t^\e],
$$
where
\begin{multline}
\bl_t^\e = \e |\p \f_t^{-\top}(x)\xi|^2 \bt_t^\e +  \label{tag}\\
+ \e \sum_{1\leq |\g| < 6} \frac{(-1)^{|\g|}}{\g!}( \p_\xi^\g |\p
\f_t^{-\top}(x)\xi|^2  )(\p_x^\g \bt_t^\e) + \e \br_6^{t,\e}.
\end{multline}
Substitution of the first term on the right hand side of \eqref{tag}
into \eqref{nextcomp} gives us precisely \eqref{nextterm}. From
\eqref{rem}, \eqref{E:unif}, and \thm{bound} we see that $\e\bFi_t
\op[\br_6^{t,\e}] = \e \bT_{t,\e}^{(1)}$, with
$\{\bT_{t,\e}^{(1)}\}$ being uniformly bounded.

Now, for all $|\g| = 1$ and any $\a$ one has
$$
| \p_x^\a( \e \p_\xi^\g |\p \f_t^{-\top}(x)\xi|^2 \, \p_x^\g
\bt_t^\e ) | \leq  C_{t,\a} \e |\xi| e^{-c_t \e |\xi|^2} \leq
\sqrt{\e} C_{t,\a}'
$$
uniformly for all $0\leq \e < \e_0$, $0\leq t <T$, $\xi \in \rn$, $x
\in \torus$. Hence, the terms with $|\g| = 1$ add up to a term of
the form $\sqrt{\e}\, \bT_{t,\e}^{(2)}$, where
$\{\bT_{t,\e}^{(2)}\}$ is uniformly bounded.

The terms with $|\g| = 2$ can be estimated as follows
$$
| \p_x^\a( \e \p_\xi^\g |\p \f_t^{-\top}(x)\xi|^2 \, \p_x^\g
\bt_t^\e ) | \leq  C_{t,\a} \e  e^{-c_t \e |\xi|^2}.
$$
So, they contribute a term $\e \bT_{t,\e}^{(3)}$. And finally, all
the terms with $|\g| > 2$ vanish.

Thus, the evolution equation \eqref{evol} has the form
\begin{equation}\label{ }
\frac{d}{dt}\bH_t^\e  = -\der{u}{\bH_t^\e } + \bA \bH_t^\e  + \e
\Delta \bH_t^\e  + \sqrt{\e}\, \bT_{t,\e}^{(4)} +
\bU_{t,\e}^{(4)},
\end{equation}
where $\{\bT_{t,\e}^{(4)}\}$ is uniformly bounded, and
$\{\bU_{t,\e}^{(4)}\}$ is uniformly compact. By the Duhamel
principle one gets
\begin{equation}\label{ }
\bH_t^\e = \bG_t^\e + \sqrt{\e}\, \int_0^t \bH_{t-s}^\e
\bT_{s,\e}^{(4)}\ ds + \int_0^t \bH_{t-s}^\e \bU_{s,\e}^{(4)}\ ds.
\end{equation}
It remains to observe that by \lem{L:unif} and \thm{bound} the
family $\{ \bH_t^\e\}$ itself is uniformly bounded, and hence, the
integrals define operators $\bT_t^\e$ and $\bU_t^\e$ with the
desired properties.

This finishes the proof of \thm{T:decomp}.

\section{Proof of the main theorem}\label{proof}
Let us recall that $\mu$ is the maximal Lyapunov exponent of the
$b$-cocycle, which determines the essential spectral radius for the
semigroup operator $\bG_t^0$ through formula \eqref{misha}. Let us
fix a $\d>0$. We can find a large $t$ such that
$$
\sup_{x,\xi} | \bB_t(x,\xi)| < \frac{1}{4} e^{t(\mu + \d)}.
$$
Then by the sharp G\"{a}rding inequality,  for $N$ large enough, we
get
\begin{equation}\label{Happr}
  \|\bH_t^\e - \bH_t^\e \bP_N\| \leq 2 \sup_{x,\xi}|\bt_t^\e(x,\xi)| <
\frac{1}{2} e^{t(\mu + \d)},
\end{equation}
for all $0 \leq \e <\e_0$. By the uniform compactness we also have
\begin{equation}\label{Uappr}
  \|\bU_t^\e - \bU_t^\e \bP_N\| <
\frac{1}{3} e^{t(\mu + \d)}.
\end{equation}

Let us fix $N$ for which both \eqref{Happr} and \eqref{Uappr} hold,
and split the semigroup $\bG_t^\e$ into the sum

\def \bGp {{\bG}^+}
\def \bGm {{\bG}^-}

\begin{equation}\label{ }
\bG_t^\e = \bGm_{t,\e} + \bGp_{t,\e},
\end{equation}
where we denote
\begin{align}
\bGm_{t,\e} & = \bH_t^\e(\bI - \bP_N) + \sqrt{\e}\, \bT_t^\e +
\bU_t^\e(\bI -  \bP_N), \\
\bGp_{t,\e} & = \bH_t^\e \bP_N + \bU_t^\e \bP_N.
\end{align}
So, $\bGp_{t,\e}$ is non-zero only on the finite-dimensional
subspace $\im \bP_N$, and in view of \eqref{Happr} and
\eqref{Uappr}, we have the estimate $\| \bGm_{t,\e} \| < \frac{5}{6}
e^{t(\mu + \d)}$ for all sufficiently small $\e$. This implies that
the resolvent $(\bGm_{t,\e} - z \bI)^{-1}$ exists and has the power
series expansion whenever $|z| > e^{t(\mu + \d)}$.

\begin{lemma}\label{L:conv}
The convergence
\begin{equation}\label{E:conv}
 \lim_{\e \ra 0}(\bGm_{t,\e} - z
 \bI)^{-1}= (\bGm_{t,0} - z \bI)^{-1}
\end{equation}
holds in the strong operator topology uniformly on compact subsets
of $\{|z|> e^{t(\mu + \d)}\}$.
\end{lemma}
\begin{proof}
Observe that $\bL^\e f \ra \bL^0 f$ for all $f \in
C^\infty(\torus)$, and $C^\infty(\torus)$ is a core of the operator
$\bL^0$. Hence, by \cite[Theorem 7.2]{Goldstein}, $\bG_t^\e \ra
\bG_t^0$ strongly. It is straightforward to prove that $\bH_t^\e \ra
\bH_t^0$ strongly, which by virtue of the decomposition
\eqref{E:decomp} also implies that $\bU_t^\e \ra \bU_t^0$. We
therefore obtain convergence $ \bGm_{t,\e} \ra \bGm_{t,0}$ in the
strong operator topology.

Thus, the conclusion of the lemma follows from the preceding
remarks.
\end{proof}

Observe that for any $|z| > e^{t(\mu +\d)}$ and $0\leq \e < \e_0$
the identity
\begin{equation}\label{ev}
    \bG_t^\e f = z f
\end{equation}
can be written as
\begin{equation}\label{ }
f + (\bGm_{t,\e} - z\bI)^{-1} \bGp_{t,\e} f = 0.
\end{equation}
This is equivalent to the system of equations
\begin{gather}
[ \bP_N + \bP_N (\bGm_{t,\e} -
z\bI)^{-1} \bGp_{t,\e} \bP_N ] f_N' = 0, \label{sys1}\\
f_N'' = -(\bGm_{t,\e} - z\bI)^{-1} \bGp_{t,\e} \bP_N f_N',
\label{sys2}
\end{gather}
where $f'_N = \bP_N f$ and $f''_N = (\bI - \bP_N) f$. We see that
$f''_N$ can be found from \eqref{sys2} if $f'_N$ is known. So, the
original eigenvalue problem \eqref{ev} is equivalent to the
finite-dimensional equation \eqref{sys1}, which in turn has a
solution at $z = z_0$ if and only if $z_0$ is a root of the analytic
function
\begin{equation}\label{ }
g(z,\e) = \det||| (e_j,e_k) +( \bP_N (\bGm_{t,\e} - z\bI)^{-1}
\bGp_{t,\e} \bP_N e_j , e_k) |||_{j,k = 1}^K,
\end{equation}
where $\{e_1,\ldots,e_K\}$ is an orthonormal basis of $\im \bP_N$.
By Lemma \ref{L:conv}, we have
\begin{equation}\label{glim}
\lim_{\e \ra 0} g(z,\e) = g(z,0)
\end{equation}
uniformly on compact sets in $\{ |z| > e^{t(\mu+\d)}\}$.

\begin{lemma}\label{L:convfull}
The resolvents $(\bG_t^\e - z \bI)^{-1}$ exist and are uniformly
bounded on compact subsets of $\{|z|> e^{t(\mu + \d)}\} \backslash
\s(\bG_t^0)$, for $0\leq \e < \e_0$, and the limit
\begin{equation}\label{E:conv}
 \lim_{\e \ra 0}(\bG_t^\e - z
 \bI)^{-1}= (\bG_t^0 - z \bI)^{-1}
\end{equation}
holds in the strong operator topology.
\end{lemma}
\begin{proof}
The existence of the resolvents follows readily from the convergence
\eqref{glim}.

Let us fix $z \not \in \s(\bG_t^0)$ and observe that
\begin{equation}\label{repr1}
    (\bG_t^\e - z \bI)^{-1} = [ \bI + (\bGm_{t,\e} - z \bI)^{-1}
    \bGp_{t,\e} ]^{-1} (\bGm_{t,\e} - z \bI)^{-1}.
\end{equation}
In view of Lemma \ref{L:conv} is suffices to show the convergence
\begin{equation}\label{conv2}
    \lim_{\e \ra 0} [ \bI + (\bGm_{t,\e} - z \bI)^{-1}
    \bGp_{t,\e} ]^{-1}  = [ \bI + (\bGm_{t,0} - z \bI)^{-1}
    \bGp_{t,0} ]^{-1}.
\end{equation}

In the direct sum $L^2 = \im \bP_N \oplus \Ker \bP_N$ we have the
following block-representation
$$
\bI + (\bGm_{t,\e} - z \bI)^{-1}  \bGp_{t,\e} = \left[
                                                  \begin{array}{cc}
                                                    \bI + \bP_N(\bGm_{t,\e} - z \bI)^{-1}
    \bGp_{t,\e} \bP_N & 0 \\
                                                    (\bI-\bP_N)(\bGm_{t,\e} - z \bI)^{-1}
    \bGp_{t,\e} \bP_N & \bI \\
                                                  \end{array}
                                                \right].
$$
So,
\begin{equation*}
[\bI + (\bGm_{t,\e} - z \bI)^{-1}  \bGp_{t,\e}]^{-1} = \left[
                                                  \begin{array}{cc}
                                                    [\bI + \bP_N(\bGm_{t,\e} - z \bI)^{-1}
    \bGp_{t,\e} \bP_N ]^{-1} & 0 \\
                            \bF_{t,\e} & \bI \\
                                                  \end{array}
                                                \right],
\end{equation*}
where
$$
 \bF_{t,\e} = - (\bI-\bP_N)(\bGm_{t,\e} - z \bI)^{-1}
    \bGp_{t,\e} \bP_N [\bI + \bP_N(\bGm_{t,\e} - z \bI)^{-1}
    \bGp_{t,\e} \bP_N ]^{-1}.
$$
Since $g(z,\e)$ is uniformly bounded away from $0$ for small $\e$,
$\bF_{t,\e}$ is uniformly bounded from above, and  hence, so is the
resolvent \eqref{repr1}. The limit \eqref{conv2} now follows from
the above formulas and \lem{L:conv}.
\end{proof}

Lemma \ref{L:convfull} already proves the spectral convergence
result for the semigroups. In order to prove it for the generators
as stated in \thm{main} we argue as follows.

Let $\l \in \s(\bL^0)$ be arbitrary. Find a $\d>0$ such that $\re \l
> \mu + \d$, and let $t>0$ be chosen as above to satisfy
\eqref{Happr} and \eqref{Uappr}. Observe the following identity:
\begin{equation}\label{ }
(\bL^\e - \zeta \bI)^{-1} = (\bG_t^\e - e^{\zeta t} \bI)^{-1}
\int_0^t e^{\zeta(t-s)} \bG_s^\e \, ds.
\end{equation}
It follows from Lemma \ref{L:convfull}, that the resolvents $(\bL^\e
- \zeta \bI)^{-1}$ are uniformly bounded on a circle $\Gamma$ of
small radius $r$ centered at $\l$ that does not contain other points
of the spectrum of $\bL^0$. Moreover,
\begin{equation}\label{conv5}
\lim_{\e \ra 0} (\bL^\e - \zeta \bI)^{-1} = (\bL^0 - \zeta
\bI)^{-1}.
\end{equation}
The Riesz projection on the spectral subspace corresponding to the
part of the spectrum of $\bL^\e$ inside $\Gamma$ is given by
\begin{equation}\label{ }
\bP^\e = \sum_{\substack{\l' \in \s(\bL^\e) \\ |\l - \l'| <r}}
\bP_{\l'}^\e =  \frac{1}{2\pi i} \int_\Gamma (\bL^\e - \zeta
\bI)^{-1} \, d \zeta.
\end{equation}
Using \eqref{conv5} the limit $\bP^\e \ra \bP_\l^0$ follows from the
dominated convergence theorem. This proves statement (iii) of our
theorem, and hence, (ii) and (i).

\subsection{Discussion}

As we noted before, the Navier-Stokes operator given by the right
hand side of \eqref{LNSE} is a particular case of the general
operator $\bL^\e$. Thus the results of \thm{main} apply to give the
convergence of the unstable eigenvalues of the Navier-Stokes
operator to eigenvalues  the Euler operator outside the essential
spectrum of the latter. Our results therefore extend the theorem of
Vishik and Friedlander \cite{VishFried} proving that a necessary
condition for instability in the Navier-Stokes equations as $\e \ra
0$ is an instability in the underlying Euler equations.

When the function space $L^2(\torus)$ is replaced by the function
space $H^m(\torus)$ it is possible to obtain in place of $\mu$ an
analogous quantity $\mu_m$ which determines the essential spectral
radius of $\bG_t^0$. The role of the $b$-cocycle defined by
\eqref{b} is replaced by a new so-called $b\xi^m$-cocycle (see
\cite{ShvCocycle,Shv2005}). All the arguments in this present paper
remain valid for the convergence of the spectrum of the viscous
operator as $\e \ra 0$ and the spectrum of the inviscid operator in
$H^m(\torus)$ with $\re \l > \mu_m$. In the particular case of the
two dimensional Euler equation in $H^1(\torus)$ it can be shown that
$\mu_1 = 0$, hence \thm{main} implies that in this example there is
precise convergence of all the points of the unstable spectra of the
Navier-Stokes operators to that of the Euler operator in the
inviscid limit.

The results of \thm{main} also apply to other fluid systems such as
the equations of geophysical fluid dynamics describing rotating,
stratified incompressible flows where the evolution operator is an
advective operator of the type $\bL^0$. We refer to \cite{Shv2005}
for an extended list of examples.

\section{Appendix}

In this section we recall a few facts about global
pseudo-differential operators (PDO) on the torus defined by
\begin{equation}\label{pdo}
  \op[\bs]f(x) = \sum_{k \in \zn} e^{i k\cdot x} \bs(x,k)
  \hat{f}(k),
\end{equation}
where $\zn = \Z^n \backslash \{0\}$, $f(x) \in \C^d$ and $\bs$ is a
$d\times d$-matrix valued symbol of class $\sym{m}$. We write $\bs
\in \sym{m}$ if $\bs \in C^\infty(\torus \times \rn)$, where $\rn =
\R^n \backslash\{0\}$, and
\begin{equation}\label{symbol}
    |\p_x^\a \p_\xi^\b \bs(x,\xi)| \leq A_{\a,\b} |\xi|^{m - |\b|},
\end{equation}
for all $|\xi| \geq 1$, $x \in \torus$, and all multi-indexes
$\a,\b$. Even though in the formula \eqref{pdo} we do not need to
require symbols to be defined outside the integer lattice, we do
assume that they are smooth in $\xi \in \rn$. For such symbols the
standard theorems of pseudo-differential calculus hold as in the
case of $\R^n$ (see \cite{DGY96}). Below we state the composition
rule with an estimate on the remainder term, which can be deduced
from a careful examination of the proof given in \cite{Wong}.

For $a,b \geq 0$, and $\{A_{\a,\b}\}$ given in \eqref{symbol}, let
us define
$$
\tilde{A}_{a,b} = \sum_{|\a| \leq a,\, |\b| \leq b} A_{\a,\b}
$$
\begin{theorem}\label{compose}
Suppose $\bs \in \sym{m_1}$ and $\bt \in \sym{m_2}$ with the
corresponding norms $\{A_{\a,\b}\}$ and $\{B_{\a,\b}\}$. Then
\begin{equation}\label{ }
\op[\bs] \circ \op[\bt] = \op[\bl],
\end{equation}
with $\bl \in \sym{m_1+m_2}$. Moreover, for all $N \in \N$, $\bl$
has the following representation
\begin{equation}\label{ }
\bl = \sum_{|\g| < N} \frac{(-1)^{|\g|}}{\g!}( \p_\xi^\g
\bs)(\p_x^\g \bt) + \br_N,
\end{equation}
where $\br_N \in \sym{m_1+m_2-N}$, and for $N>m_1+3$ satisfies the
estimate
\begin{equation}\label{rem}
    | \p_x^\a \br_N(x,\xi)| \leq c \tA{|\a|}{N+n}
    \tB{2N+|\a|-m_1-1}{0} |\xi|^{m_1+m_2+1 - N},
\end{equation}
for $|\xi| \geq 1$, where $c=c(\a,N,n,m_1,m_2)$ is independent of
the symbols.
\end{theorem}
In \eqref{rem} the restriction $N >m_1+3$ and one extra power of
$|\xi|$ is needed in order to obtain the explicit bound in terms of
the norms of $\bs$ and $\bt$. We also emphasize that the estimate
\eqref{rem} uses only the $x$-smoothness constant of $\bt$, and not
its $\xi$-smoothness.

We note that in the case of the torus the $L^2$-norm of a PDO is
bounded by the norm of only $x$-derivatives.

\begin{theorem}\label{bound}
Suppose $\bs \in \sym{0}$ satisfies \eqref{symbol}. Then $\op[\bs]$
is bounded on $L^2$ and
\begin{equation}\label{e:bound}
    \| \op[\bs] \| \leq c \tA{n+1}{0}
\end{equation}
where $c$ is independent of the symbol.
\end{theorem}
The proof uses Minkowski's inequality, and is similar to that of the
next lemma.

\begin{lemma}\label{L:uc1}
Let $\bU_\i = \op[\bs_\i]$, $\i \in I$. Suppose there exists a
constant $A >0$ independent of $\i$ such that
    \begin{equation}\label{ }
       |\p_x^\a \bs_\i(x,\xi)| \leq A |\xi|^{-1}, \quad |\xi| \geq 1,\ |\a| \leq
       n+1,
   \end{equation}
holds for all $\i \in I$. Then the family $\{ \bU_\iota\}_{\iota \in
I}$ is uniformly compact.
\end{lemma}

\begin{proof}
Let $f\in L^2$, $\|f\| = 1$, and $\supp{\hat{f}} \ss \{|k| \geq
N\}$. We obtain
\begin{align*}
\|\bU_\i f\|^2 & =\sum_{q \in \Z^n} \left| \sum_{|k|\geq N}
\hat{\bs}_\i(q-k,k) \hat{f}(k) \right|^2\\
& \lesssim \sum_{q\in \Z^n} \left| \sum_{k \not = q,\,|k|\geq N}
\hat{\bs}_\i(q-k,k) \hat{f}(k) \right|^2 + \sum_{|k| \geq N} \left|
\hat{\bs}_\i(0,k) \hat{f}(k) \right|^2 \\
& \leq A^2 \sum_{q\in\Z^n} \left( \sum_{k \not = q,\, |k|\geq N}
\frac{|\hat{f}(k)|}{|k||q-k|^{n+1}} \right)^2 + A^2 \sum_{|k|\geq N} \frac{|\hat{f}(k)|^2}{|k|^2} \\
& \leq N^{-2} A^2 \left( \sum_{q \in \zn} |q|^{-n-1} \right)^2+
N^{-2}A^2 \lesssim N^{-2}A^2.
\end{align*}
This proves the lemma.
\end{proof}

\begin{lemma}\label{L:uc2}
Let $\{\bU_\i\}_{\i\in I}$ be as in Lemma \ref{L:uc1}. Let $\bV_\k =
\op[\bt_\k]$, ${\k \in K}$, be another family such that there is a
constant $B>0$ independent of $\k$ such that
    \begin{equation}\label{tbound}
       |\p_x^\a \bt_\k(x,\xi)| \leq B , \quad |\xi| \geq 1,\ |\a| \leq
       n+1
    \end{equation}
holds for all $\k \in K$. Then the family $\{ \bU_\i \bV_\k\}_{\iota
\in I,\, \k \in K}$ is uniformly compact.
\end{lemma}

\begin{proof}
Let $N > 0$, $f\in L^2$ with $\supp \hat{f} \ss \{|k| \geq N\}$ be
fixed. Let $|q| < N/2$. Using the Cauchy-Schwartz inequality we
estimate
\begin{align*}
|\widehat{(\bV_\k f)}(q) | & = \left| \sum_{|k| \geq N}
\hat{\bt}_\k(q-k,k)\hat{f}(k) \right| \leq B  \sum_{|k| \geq N}
\frac{|\hat{f}(k)|}{|q - k|^{n+1}} \\
& \leq B \sum_{|p| > N/2} \frac{|\hat{f}(q-p)|}{|p|^{n+1}}  \leq B
\|f\| \left( \sum_{|p| > N/2} |p|^{-2(n+1)} \right)^{1/2}\\
& \lesssim N^{-1} \|f\|.
\end{align*}

So, for any fixed $M>0$ we have
$$
\lim_{N \ra \infty} \| \bP_M \bV_\k (\bI - \bP_N)\| = 0
$$
uniformly in $\k \in K$. Observe that
$$
\| \bU_\i \bV_\k - \bU_\i \bV_\k \bP_N\| \leq \| \bU_\i(\bI -
\bP_M)\bV_\k(\bI - \bP_N)\| + \| \bU_\i \bP_M \bV_\k(\bI - \bP_N)\|.
$$

Thus, using uniform compactness of $\bU_\i$'s and uniform
boundedness of $\bV_\k$'s, which follows from \thm{bound}, we can
choose $M$ large enough to make the first summand small uniformly in
$N, \i, \k$. Letting $N \ra \infty$ we make the second summand small
too. This finishes the proof.
\end{proof}


\def\cprime{$'$} \def\cprime{$'$} \def\cprime{$'$} \def\cprime{$'$}
\providecommand{\bysame}{\leavevmode\hbox
to3em{\hrulefill}\thinspace}
\providecommand{\MR}{\relax\ifhmode\unskip\space\fi MR }
\providecommand{\MRhref}[2]{%
  \href{http://www.ams.org/mathscinet-getitem?mr=#1}{#2}
} \providecommand{\href}[2]{#2}

\end{document}